\documentclass{amsart}
\usepackage{amssymb,latexsym,verbatim}
\usepackage{amsxtra,amscd}

\DeclareMathOperator{\GL}{GL}

\DeclareMathOperator{\Proj}{Proj}

\DeclareMathOperator{\Hom}{Hom}

\DeclareMathOperator{\GKdim}{GKdim}

\DeclareMathOperator{\Pic}{Pic}

\DeclareMathOperator{\QGr}{QGr}
\DeclareMathOperator{\Qch}{Qch}

\newcommand{\ov}{\overline}

\newcommand{\op}{\text{op}}

\newcommand{\Num}{\text{Num}}

\newcommand{\ZZ}{{\mathbb Z}}

\newcommand{\NN}{{\mathbb N}}
\newcommand{\PP}{{\mathbb P}}
\newcommand{\LL}{{\mathcal L}}
\newcommand{\F}{{\mathcal F}}

\newcommand{\calC}{{\mathcal C}}

\newcommand{\calH}{{\mathcal H}}

\newcommand{\calM}{{\mathcal M}}

\newcommand{\OO}{{\mathcal O}}

\newcommand{\gs}{{\sigma}}

\numberwithin{equation}{section}

\newtheorem{theorem}{Theorem}[section]
\newtheorem{corollary}[theorem]{Corollary}
\newtheorem{lemma}[theorem]{Lemma}
\newtheorem{proposition}[theorem]{Proposition}

\theoremstyle{definition}

\newtheorem{definition}[theorem]{Definition}

\newtheorem{hypothesis}[theorem]{Hypothesis}

\theoremstyle{remark}

\begin{document}

\setcounter{page}{1}

\title{Noncommutative ampleness for multiple divisors}

\author{Dennis S. Keeler}
     \thanks{ %Preprint. To be submitted for publication.
     Partially supported by NSF grant DMS-9801148 and NSF Postdoctoral Fellowship.}  
\address{ Department of Mathematics \\ MIT \\ Cambridge, MA 02139-4307   }
\email{dskeeler@mit.edu}
\urladdr{http://www.mit.edu/\~{}dskeeler}
 
%\date{ \today }


% Primary
% 14A Foundations
%	14A22 Noncommutative algebraic geometry
% 14C Cycles and subschemes
%	14C20 Divisors, linear systems, invertible sheaves
% 14F (Co)homology theory
%	14F05 Vector bundles, sheaves, related constructions
%	14F17 Vanishing theorems
% 14J Surfaces and higher-dimensional varieties
%	14J50 Automorphisms of surfaces and higher-dimensional varieties
% 16P Chain conditions, growth conditions, and other forms of finiteness
%	16P40 noetherian rings and modules
%	16P90 Growth rate, Gelfand-Kirillov dimension
% 16S Rings and algebras arising under various constructions
%	16S38 Rings arising from non-commutative algebraic geometry
%	16S80 Deformations of rings
% 16W Rings and algebras with additional structure
%	16W50 Graded rings and modules
%
% Secondary
% 14C Cycles and subschemes
%	14C17 Intersection theory, characteristic classes, intersection multiplicities
%	14C40 Riemann-Roch theorems
%	

\begin{abstract}
The twisted homogeneous coordinate ring is one of the basic constructions of the
noncommutative projective geometry of Artin, Van den Bergh, and others. Chan
generalized this construction to the multi-homogeneous case, using a concept of
right ampleness for a finite collection of invertible sheaves and automorphisms
of a projective scheme. From this he derives that certain multi-homogeneous rings,
such as tensor products of twisted homogeneous coordinate rings, are right noetherian.
We show that right and left ampleness are equivalent and that there is a simple
criterion for such ampleness. Thus we find under natural hypotheses that 
multi-homogeneous coordinate rings are noetherian and have integer GK-dimension.
\end{abstract}

\maketitle

\section{Introduction}

Let $R$ be an $\NN$-graded algebra over an algebraically closed field $k$
such that $\dim R_i < \infty$ for all $i$.
One of the main techniques of noncommutative projective geometry is to study
a graded ring $R$ via a category $\calC$ of graded $R$-modules.
More specifically, one usually examines $\QGr R$, the quotient category of graded
right
$R$-modules modulo the full subcategory of torsion modules; one hopes
that $\QGr R$ will have geometric properties, since the Serre Correspondence
Theorem says that if $R$ is commutative and generated in degree one, then
there is a category equivalence
$\QGr R \cong \Qch X$, where $\Qch X$ is the category of quasi-coherent
sheaves on $X = \Proj R$ \cite[Exercise~II.5.9]{Ha1}.

The twisted homogeneous coordinate rings are the 
most basic class of rings in noncommutative geometry.
Such a ring $R$ is constructed from a commutative projective scheme $X$,
an automorphism $\gs$ of $X$, and an invertible sheaf $\LL$.
When the pair $(\LL, \gs)$ satisfies ``right $\sigma$-ampleness,''
then $R$ is right noetherian and
has $\QGr R \cong \Qch X$ \cite[Theorems~1.3, 1.4]{AV}.
These rings were first used to show that Artin-Schelter regular algebras
of dimension $3$ are noetherian domains \cite{ATV,St1,St2} and their basic properties
were studied in \cite{AV}. Further, any domain of GK-dimension $2$, generated
in degree one, is a twisted homogeneous coordinate ring for some curve $X$
\cite{AS}.

A simple criterion for right $\sigma$-ampleness was found in \cite{Ke1}.
From this criterion one sees that right and left $\gs$-ampleness
are equivalent. Hence the associated ring $R$ is noetherian. One also
sees that the GK-dimension of $R$ is an integer. (While this paper and
\cite{Ke1} work over an algebraically closed field, we note
that \cite{Ke-filters} generalized these results to the case
of a commutative noetherian base ring.)

Chan introduced twisted multi-homogeneous coordinate rings in \cite{Chan},
which are constructed from a finite collection $\{ (\LL_i, \gs_i) \}$
of invertible sheaves and automorphisms
 on a projective scheme $X$.
When the set $\{ (\LL_i, \gs_i) \}$ is ``right ample,'' 
then the category $\QGr R$ of multi-graded right $R$-modules
modulo torsion modules again has $\QGr R \cong \Qch X$.
With some natural extra hypotheses, $R$ will be right noetherian.
Via these methods, Chan shows that some rings associated to twisted
homogeneous coordinate rings, like tensor products of two such coordinate
rings, are right noetherian.

In this paper, we will generalize the results of \cite{Ke1} to the multi-homogeneous
case and thereby strengthen \cite{Chan}.
More specifically, we show

\begin{theorem}[see Theorem~\ref{th:multi-main}, Corollary~\ref{cor:multi-leftright}]
\label{th:main-intro}
Let $X$ be a projective scheme and let $\{ (\LL_i, \gs_i) \}$ be a finite set of pairs
of invertible sheaves and automorphisms. Then there is a simple criterion
for $\{ (\LL_i, \gs_i) \}$ to be right ample.
This criterion shows that right and left ampleness are equivalent.
\end{theorem}

We then immediately have, in Corollary~\ref{cor:tensor-product}, that
the tensor product $B \otimes_k B'$
is noetherian, where $B,B'$ are twisted homogeneous coordinate rings associated
to ample pairs $(\LL, \gs),(\LL', \gs')$. If $B$ is generated in degree one
 and $I$ is the irrelevant ideal
of $B$, then the Rees algebra $B[I t]$ is noetherian; see Corollary~\ref{cor:ReesRing}.

We also show

\begin{theorem}[see Theorem~\ref{th:GK}]  \label{th:GK-intro}
Let $B$ be a twisted multi-homogenous coordinate ring
under suitable hypotheses~\eqref{hyp:ringB}. Then
$\GKdim B$ is an integer with geometrically defined bounds.
\end{theorem}

Most of this paper appeared in the author's Ph.D. thesis, under the direction
of J.T.~Stafford.

\section{Right ampleness is left ampleness}

Because of the notational difficulties associated with
handling the ampleness of arbitrarily many pairs $(\LL_i, \gs_i)$, we will use the concept  
of an invertible bimodule $\LL_\gs$. % in Definition~\ref{def:invertible-bimodule}.
In this paper it will only be important to know how
invertible bimodules act on a coherent sheaf $\F$, so we will treat $\LL_\gs$
as a notational convenience where
\[
\F \otimes \LL_\gs = \gs_*(\F \otimes \LL),\quad \LL_\gs \otimes \F = \LL \otimes \gs^*\F
\]
and the right hand side of the above equations are 
just $\OO_X$-modules. For a formal definition of invertible bimodule see
\cite[$\S$2]{AV}.
\begin{comment}
Let $\gs$ and $\tau$ be automorphisms
and let $\LL$ and $\calM$ be invertible sheaves. Then for a coherent sheaf $\F$,
\begin{align*}
\LL_\gs \otimes (\calM_\tau \otimes \F) &= \LL_\gs \otimes (\calM \otimes \tau^*\F) \\
 &= \LL \otimes \gs^*\calM \otimes \gs^* \tau^* \F \\
 &= \LL \otimes \gs^*\calM \otimes (\tau \gs)^*\F \\
 &= (\LL \otimes \gs^*\calM)_{\tau \gs} \otimes \F.
\end{align*}
\end{comment}
Given
two invertible bimodules $\LL_\gs$ and $\calM_\tau$, one finds the tensor product to
be
\begin{equation}\label{eq:bimodproduct}
\LL_\gs \otimes \calM_\tau = (\LL \otimes \gs^* \calM)_{\tau \gs},
\end{equation}
where the second tensor product is the usual product on quasi-coherent sheaves
\cite[Lemma~2.14]{AV}.
We will sometimes denote the product of invertible bimodules by
juxtaposition if the meaning is clear.
The $\OO_X$-module underlying a product of bimodules $\LL_\gs \otimes \calM_\tau$
will be denoted $\vert \LL_\gs \otimes \calM_\tau \vert$; in this particular
case $\vert \LL_\gs \otimes \calM_\tau \vert = \LL \otimes \gs^*\calM$.

We will also use the notation $\LL^\gs = \gs^*\LL$. The automorphism $\gs$ induces
a natural isomorphism
\begin{equation}\label{eq:switch sides}
\F \otimes \LL_\gs = \gs_*(\F \otimes \LL) \cong \LL^{\gs^{-1}} \otimes \F^{\gs^{-1}} = 
\LL^{\gs^{-1}}_{\gs^{-1}} \otimes \F,
\end{equation}
for any coherent sheaf $\F$.

We now sketch the construction of a twisted multi-homogeneous coordinate ring;
for details see \cite[$\S$2]{Chan}.
Let $\{ (\LL_i)_{\gs_i} \}$ be a collection of $s$ invertible bimodules,
possibly with repetitions.
For notational convenience, we will write $\LL_{(i,\gs_i)} = (\LL_i)_{\gs_i}$.
Given these $s$ invertible bimodules, 
one wishes to form an associated twisted multi-homogeneous
coordinate ring $B = B(X; \{ \LL_{(i, \gs_i)} \} )$. 
For an $s$-tuple $\ov{n}= (n_1, \dots, n_s)$ we define the multi-graded piece $B_{\ov{n}}$
as
\begin{equation}\label{eq:multi-pieces}
B_{\ov{n}} = H^0 (X, \LL_{(1,\gs_1)}^{n_1} \dots \LL_{(s,\gs_s)}^{n_s})
\end{equation}
where the cohomology of an invertible bimodule is just cohomology of the underlying sheaf.
Multiplication should be given by
\begin{equation}\label{eq:multi-multiplication}
a \cdot b = a {{\gs}^{\ov{m}}}(b)
\end{equation}
when $a \in B_{\ov{m}}$ and $b \in B_{\ov{n}}$.
Here ${{\gs}^{\ov{m}}}(b) = \gs_1^{m_1} \gs_2^{m_2} \dots \gs_s^{m_s} (b)$,
where the action of an automorphism on a global section is induced by pullback.

\begin{comment}
However, for this multiplication to be defined, the invertible bimodules must commute with
each other. To see this, consider the bigraded case, with bimodules $\LL_\gs, \calM_\tau$.
Then 
\[
B_{(1,0)} = H^0(X, \LL_\gs),\quad B_{(0,1)} = H^0(X, \calM_\tau),\quad 
B_{(1,1)} = H^0(X, \LL_\gs \calM_\tau).
\]
Given the multiplication above, we have $B_{(1,0)} B_{(0,1)} \subset B_{(1,1)}$
and $B_{(0,1)} B_{(1,0)} \subset H^0(X, \calM_\tau \LL_\gs)$.
To guarantee that $B_{(1,1)} = H^0(X, \calM_\tau \LL_\gs)$, so that we
have a bigraded ring, we demand $\LL_\gs \calM_\tau = \calM_\tau \LL_\gs$.
\end{comment}

However, to make the ring construction work, \cite{Chan} shows that we need
the invertible bimodules to commute with each other.
 Examining \eqref{eq:bimodproduct}, we see that two bimodules $\LL_\gs, \calM_\tau$ 
commute when
\begin{equation}\label{eq:bimod-commute}
\LL \otimes \gs^* \calM \cong \calM \otimes \tau^* \LL\quad \text{and} \quad \gs\tau = \tau\gs.
\end{equation}
Thus we need sheaf isomorphisms 
$\varphi_{ij}: \LL_{(j, \gs_j)} \LL_{(i, \gs_i)} \to \LL_{(i, \gs_i)} \LL_{(j,\gs_j)}$ for
each $1 \leq i < j \leq s$.
It is further noted in \cite{Chan} that when there are three or more bimodules, these 
isomorphisms must be compatible on ``overlaps'' in the sense of Bergman's Diamond Lemma.
In terms of the isomorphism $\varphi_{ij}$ this means \cite[p.~444]{Chan}
\begin{multline}\label{eq:Bergman}
(\varphi_{ij} \otimes 1_{\LL_{(k, \gs_k)}}) 
	\circ (1_{\LL_{(j, \gs_j)}} \otimes \varphi_{ik})
	\circ (\varphi_{jk} \otimes 1_{\LL_{(i, \gs_i)}} )\\
= ( 1_{\LL_{(i, \gs_i)}} \otimes \varphi_{jk}) 
	\circ (\varphi_{ik} \otimes 1_{\LL_{(j, \gs_j)}})
	\circ (1_{\LL_{(k, \gs_k)}} \otimes \varphi_{ij})
\end{multline}
in $\Hom(\LL_{(k, \gs_k)} \LL_{(j, \gs_j)} \LL_{(i, \gs_i)}, 
\LL_{(i, \gs_i)} \LL_{(j, \gs_j)} \LL_{(k, \gs_k)})$. 
We will always assume that we have this compatibility when forming the ring $B$.
Summarizing, we have
\begin{proposition}\label{prop:multi-def}
Let $\{ \LL_{(i, \gs_i)} \}$ be a finite collection of commuting invertible bimodules. Assume that
these bimodules
have compatible pairwise commutation relations in the sense of \eqref{eq:Bergman}.
Then there is a multi-graded ring $B$ with multi-graded pieces given
by \eqref{eq:multi-pieces} and multiplication given by \eqref{eq:multi-multiplication}.\qed
\end{proposition}

To study these rings, a multi-graded version of $\gs$-ampleness is introduced.
Since we will be interested in both this version of ampleness and the usual commutative
one, we will
call this (right) \emph{NC-ampleness}, 
whereas \cite{Chan} uses the terminology (right) ampleness. 
We define the  ordering on $s$-tuples
to be the standard one, i.e., $(n_1',\dots, n_s') \geq (n_1, \dots, n_s)$ if $n_i' \geq n_i$
for all $i$. 
For simplicity we write 
$\LL_{\ov{\gs}}^{\ov{m}} = \LL^{m_1}_{(1, \gs_1)} \dots \LL^{m_s}_{(s, \gs_s)}$.

\begin{definition}\label{def:NC-ample}
Let $X$ be a projective scheme with $s$ commuting invertible 
bimodules $\{ \LL_{(i, \gs_i)} \}$. 
\begin{enumerate}
\item If for any coherent sheaf $\F$, there exists an $\ov{m}_0$ 
	such that
		\[
		H^q(X, \F \otimes \LL_{\ov{\gs}}^{\ov{m}}) = 0
		\] for  $q > 0$ and $\ov{m} \geq \ov{m}_0$, then the
		set $\{ \LL_{(i, \gs_i)} \}$ is called \emph{right NC-ample}.
\item If for any coherent sheaf $\F$, there exists an $\ov{m}_0$ 
	such that
		\[
		H^q(X, \LL_{\ov{\gs}}^{\ov{m}} \otimes \F ) = 0
		\] for  $q > 0$ and $\ov{m} \geq \ov{m}_0$, then the
		set $\{ \LL_{(i, \gs_i)} \}$ is called \emph{left NC-ample}.
\end{enumerate}
\end{definition}

As in the case of one invertible bimodule, right and left NC-ampleness are
related.

\begin{lemma}\label{lem:multi-rightleft}
{\upshape (cf. \cite[Lemma~2.3]{Ke1})} %Lemma~\ref{lem:rightleft} 
Let $X$ be a projective scheme with $s$ commuting invertible 
bimodules $\{ (\LL_i)_{\gs_i} \}$.  Then the
set $\{ (\LL_i^{\gs_i^{-1}})_{\gs_i^{-1}} \}$ commutes pairwise.
Also, the set $\{ (\LL_i)_{\gs_i} \}$ is right NC-ample if and only if
the set $\{ (\LL_i^{\gs_i^{-1}})_{\gs_i^{-1}} \}$ is left NC-ample.
\end{lemma}
\begin{proof}
Let $\LL_\gs, \calM_\tau$ be two commuting invertible bimodules.
Then \eqref{eq:bimod-commute} holds. Obviously $\gs^{-1} \tau^{-1} = \tau^{-1} \gs^{-1}$.
Now since $\LL \otimes \gs^* \calM \cong \calM \otimes \tau^* \LL$, pulling back
by $(\gs^{-1} \tau^{-1})$ we have
\[
(\tau^{-1})^* (\gs^{-1})^* \LL \otimes (\tau^{-1})^* \calM 
\cong (\gs^{-1})^* (\tau^{-1})^* \calM \otimes (\gs^{-1})^* \LL.
\]
So $\LL^{\gs^{-1}}_{\gs^{-1}} = ((\gs^{-1})^* \LL)_{\gs^{-1}}$ and  
$\calM^{\tau^{-1}}_{\tau^{-1}} = ((\tau^{-1})^* \calM)_{\tau^{-1}}$
commute.

Now using \eqref{eq:switch sides} and the fact that the bimodules commute, 
we see that  
\[
H^q(X, \F \otimes  (\LL_1)^{m_1}_{\gs_1}\dots (\LL_s)^{m_s}_{\gs_s}) =
H^q(X, (\LL_1^{\gs_1^{-1}})^{m_1}_{\gs_1^{-1}}\dots (\LL_s^{\gs_s^{-1}})^{m_s}_{\gs_s^{-1}}
\otimes \F )
\]
for all $q, m_i$. Thus right NC-ampleness of $\{ (\LL_i)_{\gs_i} \}$  is equivalent
to left NC-ampleness of $\{ (\LL_i^{\gs_i^{-1}})_{\gs_i^{-1}} \}$.
\end{proof}

%We also have an analogue of Lemma~\ref{lem:Bgs-Bgsinv-op}.

\begin{lemma}\label{lem:opposite}
Let $X$ be a projective scheme over $k$ with $s$ commuting invertible 
bimodules $\{ (\LL_i)_{\gs_i} \}$. Assume that the commutation relations
of $\{ (\LL_i)_{\gs_i} \}$ and of $\{ (\LL_i^{\gs_i^{-1}})_{\gs_i^{-1}} \}$ 
are compatible in the sense of \eqref{eq:Bergman}. 
If $B' = B(X; \{ (\LL_i)_{\gs_i} \})$ and 
$B= B(X; \{ (\LL_i^{\gs_i^{-1}})_{\gs_i^{-1}} \})$, then
$B \cong (B')^\op$.
\end{lemma}
\begin{proof}
Let
$\tau\colon B \to (B')^\op$ given by $\tau(a) = \gs_1^{n_1}\dots\gs_s^{n_s}(a)$
for $a \in B_{(n_1, \dots, n_s)}$.
Extend $\tau$ linearly so it is a vector space map. It is obviously a
vector space isomorphism. 

Let $\cdot$ be multiplication in $B$ and $\ast$ be multiplication in $(B')^\op$.
For $a \in B_{\ov{n}}, b \in B_{\ov{m}}$,
\begin{align*}
\tau(a \cdot b) &= \tau(a\gs^{-\ov{n}}(b)) = \gs^{\ov{n}+\ov{m}}(a) \gs^{\ov{m}}(b), \\
\tau(a) \ast \tau(b) &= \gs^{\ov{n}}(a) \ast \gs^{\ov{m}}(b) = \gs^{\ov{m}}(b) \gs^{\ov{n}+\ov{m}}(a).
\end{align*}
Thus $\tau(a \cdot b) = \tau(a) \ast \tau(b)$, as required.
\end{proof}

As in \cite[Proposition~2.3]{Ke1}, we have simpler equivalent conditions
for a set of bimodules to be right NC-ample.

\begin{proposition} \label{prop:multi-tfae}
Let $X$ be a projective scheme with $s$ commuting invertible 
bimodules $\{ \LL_{(i, \gs_i)} \}$. Then
the following are equivalent:
\begin{enumerate}
	\item \label{prop:multi-tfae1} The set $\{ \LL_{(i, \gs_i)} \}$
	is right NC-ample.
	\item \label{prop:multi-tfae2} For any coherent sheaf $\F$, 
	there exists an $\ov{m}_0$ such that
		$\F \otimes \LL_{\ov{\gs}}^{\ov{m}}$ is 
		generated by global sections for $\ov{m} \geq \ov{m}_0$.
	\item \label{prop:multi-tfae3} For any invertible sheaf $\calH$, 
	there exists an $\ov{m}_0$ such that
		$\vert \calH^{-1} \otimes \LL_{\ov{\gs}}^{\ov{m}} \vert$ 
		is very ample for $\ov{m} \geq \ov{m}_0$.
	\item \label{prop:multi-tfae4} For any invertible sheaf $\calH$, 
	there exists an $\ov{m}_0$ such that
		$\vert \calH^{-1} \otimes \LL_{\ov{\gs}}^{\ov{m}} \vert$ 
		is  ample for $\ov{m} \geq \ov{m}_0$.
\end{enumerate}
A similar statement holds for left NC-ample.
\end{proposition}

\begin{proof}
This is a special case of \cite[Theorem~1.3, Proposition~6.9]{Ke-filters}.
\end{proof}

We can now give a connection between right NC-ampleness and the
concept of $\gs$-ampleness for one invertible sheaf $\LL$.

\begin{lemma}\label{lem:NC-ample-then-s-ample}
Let $X$ be a projective scheme with $s$ commuting invertible 
bimodules $\{ \LL_{(i, \gs_i)} \}$. Suppose that $\ov{n} = (n_1, \dots, n_s) \in (\NN^+)^s$
and set $\tau = \gs_1^{n_1}\dots\gs_s^{n_s}$.
If the set of bimodules is right NC-ample,
then $\vert \LL_{(1, \gs_1)}^{n_1}\dots \LL_{(s, \gs_s)}^{n_s} \vert$
is $\tau$-ample.
\end{lemma}
\begin{proof}
Let $\calH$ be an invertible sheaf and
let $\ov{m}_0$ be such that for all $\ov{m} \geq \ov{m}_0$,
the sheaf $\vert \calH^{-1} \otimes \LL_{\ov{\gs}}^{\ov{m}} \vert$ is ample
by Proposition~\ref{prop:multi-tfae}\eqref{prop:multi-tfae4}.

Now there exists an integer $l_0$ such that for all $l \geq l_0$,
we have $l \ov{n} \geq \ov{m}_0$.
So $\vert \calH^{-1} \otimes (\LL_{\ov{\gs}}^{\ov{n}})^l \vert$
is ample. Thus by \cite[Proposition~2.3(4)]{Ke1},
$\vert \LL_{\ov{\gs}}^{\ov{n}} \vert$ is $\tau$-ample.
\end{proof}

Recall that the Picard group of $X$ modulo numerical equivalence, $A^1_{\Num}(X)=\Pic X/\equiv$,
is a finitely generated free abelian group \cite[p.~305, Remark 3]{Kleiman}.
Thus the action of $\gs$ on $A^1_{\Num}(X)$ is given by some $P \in \GL_\rho(\ZZ)$
for some $\rho > 0$.
We say that $\gs$ is \emph{unipotent} if all the eigenvalues of $P$ equal $1$
and that $\gs$ is \emph{quasi-unipotent} if all the eigenvalues of $P$ are roots of unity.
This is a well-defined notion \cite[Proposition~7.12]{Ke-filters}.
We then  have a new version
of \cite[Theorem~1.3]{Ke1}.

\begin{theorem}\label{th:multi-main}
Let $X$ be a projective scheme with $s$ commuting invertible bimodules $\{ \LL_{(i,\gs_i)} \}$.
The set $\{ \LL_{(i,\gs_i)} \}$ is (right) NC-ample if and only if every $\gs_i$ is
quasi-unipotent and there exists $\ov{m}_0 \in \NN^s$ such that
$\vert \LL_{\ov{\gs}}^{\ov{m}} \vert$ is ample for  all $\ov{m} \geq \ov{m}_0$.
\end{theorem}

\begin{proof}
Suppose that $\{ \LL_{(i,\gs_i)} \}$ is right NC-ample. Then by 
Proposition~\ref{prop:multi-tfae}\eqref{prop:multi-tfae4},
there exists $\ov{m}_0 \in \NN^s$ such that
$\vert \LL_{\ov{\gs}}^{\ov{m}} \vert$ is ample for  all $\ov{m} \geq \ov{m}_0$.
Further, by the previous lemma, 
$\LL^{n_1}_{(1,\gs_1)}\dots \LL^{n_s}_{(s,\gs_s)}$ is $\tau$-ample when
$\tau = \gs_1^{n_1} \dots \gs_s^{n_s}$ and each $n_i > 0$.
%Then one has that $\tau$ is quasi-unipotent by Theorem~\ref{th:main}.
Now recall that all the automorphisms commute and hence their actions on $A^1_{\Num}(X)$
are commuting matrices. Thus the eigenvalues of the product $\gs_1^{n_1} \dots \gs_s^{n_s}$
 are products of eigenvalues from each $\gs_i$.
So if $\gs_1$ were not quasi-unipotent, then either $\tau_1 = \gs_1 \gs_2 \dots \gs_s$ or
$\tau_2 = \gs_1^2 \gs_2 \dots \gs_s$ would not be quasi-unipotent.
But $\tau_1$ and $\tau_2$ must be quasi-unipotent by \cite[Theorem~1.3]{Ke1}
since the corresponding sheaves $\LL^{1}_{(1,\gs_1)}\dots \LL^{1}_{(s,\gs_s)}$
and $\LL^{2}_{(1,\gs_1)}\dots \LL^{1}_{(s,\gs_s)}$ are $\tau_1$-ample and $\tau_2$-ample
respectively. Thus each $\gs_i$ must be quasi-unipotent.

Now suppose that every $\gs_i$ is
quasi-unipotent and there exists $\ov{m}_0 \in \NN^s$ such that
$\vert \LL_{\ov{\gs}}^{\ov{m}} \vert$ is ample for  all $\ov{m} \geq \ov{m}_0$.
As the $\gs_i$ commute, $\tau = \gs_1 \dots \gs_s$ is quasi-unipotent.
Then by \cite[Theorem~1.3]{Ke1}, the invertible bimodule $\LL_{(1,\gs_1)} \dots \LL_{(s,\gs_s)}$
is $\tau$-ample. So given any invertible sheaf
$\calH$, there exists $n_0 \in \NN$ such that 
\[
\vert \calH^{-1} \otimes (\LL_{(1,\gs_1)} \dots \LL_{(s,\gs_s)})^n \vert =
\vert \calH^{-1} \otimes \LL_{(1,\gs_1)}^n \dots \LL_{(s,\gs_s)}^n \vert
\]
is ample for $n \geq n_0$ by \cite[Proposition~2.3(4)]{Ke1}.
Then for all $\ov{m} \geq (n_0, n_0, \dots, n_0) + \ov{m}_0$ the invertible
sheaf
\[
\vert \calH^{-1} \otimes \LL_{\ov{\gs}}^{\ov{m}} \vert = 
\vert \calH^{-1} \otimes \LL_{(1,\gs_1)}^{n_0} \dots \LL_{(s,\gs_s)}^{n_0} \vert \otimes
\vert \LL_{(1,\gs_1)}^{m_1 - n_0} \dots \LL_{(s,\gs_s)}^{m_s - n_0} 
\vert^{\gs_1^{n_0} \dots \gs_s^{n_0}}
\]
is the tensor product of two ample invertible sheaves.
Hence it
is ample and so
the set of invertible bimodules is right NC-ample 
by Proposition~\ref{prop:multi-tfae}\eqref{prop:multi-tfae4}.
\end{proof}

\begin{corollary}\label{cor:multi-leftright}
Let $X$ be a projective scheme with $s$ commuting invertible bimodules
$\{ \LL_{(i, \gs_i)} \}$. Then  $\{ \LL_{(i, \gs_i)} \}$ is right NC-ample if and only if
it is left NC-ample.
\end{corollary}
\begin{proof}
Suppose that $\{ \LL_{(i, \gs_i)} \}$ is right NC-ample. 
Then each $\gs_i$ is quasi-unipotent and
there exists $\ov{m}_0$ such that 
$\vert \LL_{(1, \gs_1)}^{m_1}\dots\LL_{(s, \gs_s)}^{m_s} \vert$ 
is ample for $(m_1, \dots, m_s) \geq \ov{m}_0$.
Pulling back by $\gs_1^{-m_1}\dots \gs_s^{-m_s}$, we have that
$\vert (\LL_1^{\gs_1^{-1}})_{\gs_1^{-1}}^{m_1}
 \dots (\LL_s^{\gs_s^{-1}})_{\gs_s^{-1}}^{m_s} \vert$
 is ample. Thus by Theorem~\ref{th:multi-main}, the set
 $\{ (\LL_i^{\gs_i^{-1}})_{\gs_i^{-1}} \}$ is right NC-ample.
So the original set $\{ \LL_{(i, \gs_i)} \}$ is left NC-ample
by Lemma~\ref{lem:multi-rightleft}.
The argument is clearly reversible.
\end{proof}

Thus we may now refer to a set of bimodules as being simply NC-ample.

Note the difference between \cite[Theorem~1.3]{Ke1} and Theorem~\ref{th:multi-main}.
The former requires only that $\vert \LL_\gs^m \vert$ is ample for one value
of $m$, while the latter requires the product of bimodules to be ample for all $\ov{m} \geq \ov{m}_0$.
To see this stronger requirement is necessary, let $X$ be any projective
scheme with $\LL$ any ample invertible sheaf.
We need to rule out  the pair $\LL, \LL^{-1}$ where
 the bimodule action is
the usual commutative one. In this particular case, of course 
$\LL^1 \otimes (\LL^{-1})^0$ is ample. But $\LL^{m_1} \otimes (\LL^{-1})^{m_2}$
is not ample for all $(m_1, m_2)$ sufficiently large; just fix $m_1$
and let $m_2$ go to infinity.

It is not necessary for one of the $\LL_{(i, \gs_i)}^m$ to 
be ample for $m \gg 0$, since on $\PP^1 \times \PP^1$, the pair $\OO(1,0), \OO(0,1)$
is NC-ample, where again these bimodules act only as commutative
invertible sheaves.

\section{Ring theoretic consequences}

Unlike the case of only one bimodule, the multi-graded ring $B$ may not be noetherian
when $\{ \LL_{(i,\gs_i)} \}$ is NC-ample. In fact, \cite[Example~5.1]{Chan} gives
a simple commutative (and hence not finitely generated) counterexample.
However, Chan
introduces an additional property for an invertible bimodule $\LL_\gs$ on
$X$ to guarantee the noetherian condition.

\begin{hypothesis}\label{property*} 
There exists a projective scheme $Y$ with automorphism $\gs$ and a $\gs$-equivariant
morphism $f\colon X \to Y$. That is $\gs_Y \circ f = f \circ \gs_X$. There also
exists an invertible sheaf $\LL'$ on $Y$ such that $\LL = f^* \LL'$ and such
that $\LL'_\gs$ is $\gs$-ample. (\cite{Chan} labels this property $(\ast)$.)
\end{hypothesis}

This property \eqref{property*} is saying that for $m \gg 0$,
$\vert \LL_\gs^m \vert$ is generated by global sections, since it is 
a pullback of $\vert (\LL')^m_\gs \vert$, which is eventually very ample
by \cite[Proposition~2.3(3)]{Ke1}.
Note in particular that if $\LL$ is already $\gs$-ample, then
$\LL_\gs$ satisfies \eqref{property*} trivially.
 Using this property,
one determines

\begin{theorem}\label{th:ChanNoetherian}
{\upshape (\cite[Theorem~5.2]{Chan})}
Let $X$ be a projective scheme with commuting invertible bimodules $\LL_\gs, \calM_\tau$.
Suppose that the pair is (right) NC-ample and each bimodule satisfies \eqref{property*}, 
possibly for different $Y$.
Then $B(X; \LL_\gs, \calM_\tau)$ is right noetherian.\qed
\end{theorem}

Then combining Corollary~\ref{cor:multi-leftright}, Lemma~\ref{lem:opposite},
 and the theorem above, we have
\begin{theorem}\label{th:multi-noeth}
Let $X$ be a projective scheme with commuting invertible bimodules $\LL_\gs, \calM_\tau$.
Suppose that the pair is NC-ample and each bimodule satisfies \eqref{property*},
possibly for different $Y$.
Then $B(X; \LL_\gs, \calM_\tau)$ is noetherian.\qed
\end{theorem}

Now we can prove that two particularly interesting twisted multi-homogeneous 
coordinate rings, a Rees ring and a tensor product, are noetherian,
strengthening the results of \cite[Corollaries~5.7, 5.8]{Chan}.
In the latter case, we may replace his proof, based on spectral sequences, 
by an easier one
since the criterion of Theorem~\ref{th:multi-main}
simplifies testing the NC-ampleness of the relevant pair of bimodules.

\begin{corollary}\label{cor:ReesRing} %{\upshape (\cite[Corollary~5.7]{Chan})}
Let $\LL_\gs$ be a $\gs$-ample invertible bimodule on a projective scheme $X$. 
Let $B = B(X; \LL_\gs)$ be generated in degree one.
Then
the Rees algebra $B[It] = \oplus_{r=0}^\infty I^r t^r$ of $B$ is noetherian,
where $I= B_{> 0}$ is the irrelevant ideal.
\end{corollary}
\begin{proof}
The ring $B[It]$ has bigraded pieces
\[
B_{(i, j)} = H^0(X, \LL_\gs^i \LL_\gs^j) t^j
\]
since $I^j = \oplus_{l = j}^\infty B_l$ when $B$ is generated in degree one.
The pair $\LL_\gs, \LL_\gs$ is obviously NC-ample and satisfies \eqref{property*}.
Thus Theorem~\ref{th:multi-noeth} applies.
\end{proof}

%In the following corollary, we not only get a stronger result of
%noetherian versus just right noetherian, but it is also
%easier to test the NC-ampleness of the relevant pair of bimodules.

\begin{corollary}\label{cor:tensor-product} %{\upshape (\cite[Corollary~5.8]{Chan})}
Let $\LL_\gs$ be $\gs$-ample on a projective scheme $X$ and let $\calM_\tau$ be
$\tau$-ample on a projective scheme $Y$. Then
$B(X; \LL_\gs) \otimes B(Y; \calM_\tau)$ is noetherian.
\end{corollary}

\begin{proof}
It is argued in \cite[Example~4.3]{Chan} that
\[
B(X; \LL_\gs) \otimes B(Y; \calM_\tau) \cong
B(X \times Y; (\pi_1^* \LL)_{\gs \times 1}, (\pi_2^* \calM)_{1 \times \tau}),
\]
where the $\pi_i$ are the natural projections.
These two invertible bimodules on $X \times Y$ obviously satisfy \eqref{property*}.
%It suffices to show the pair is NC-ample, which is shown in
%\cite[Corollary~5.8]{Chan} using spectral sequences. However, we can now given
%a simpler proof using Theorem~\ref{th:multi-main}.

Since $\LL_\gs$ is $\gs$-ample
and $\calM_\tau$ is $\tau$-ample, there is an $m_0$ such that
$\vert \LL_\gs^m \vert$ and $\vert \calM_\tau^m \vert$ is ample for all $m \geq m_0$.
Note that $(\gs \times 1)^* \pi_1^* \LL = \pi_1^* \gs^* \LL$ and a similar
formula holds for $\calM_\tau$. Then
\[
\vert (\pi_1^* \LL)_{\gs \times 1}^{m_1} (\pi_2^* \calM)_{1 \times \tau}^{m_2} \vert
\]
is ample for all $(m_1, m_2) \geq (m_0, m_0)$ by
\cite[p.~125, Exercise~5.11]{Ha1}.

Now $\gs$ is quasi-unipotent and we wish to show $\gs \times 1$ is as well.
It is tempting to think that as a matrix acting on $A^1_{\Num}(X \times Y)$ one has
$\gs \times 1 = \gs \oplus 1$. However, this may not be the case, since in
general $A^1_{\Num}(X \times Y)$ has larger rank than $A^1_{\Num}(X) \oplus A^1_{\Num}(Y)$
\cite[p.~367, Exercise~1.6]{Ha1}.
But let $\calH_X$ and $\calH_Y$ be ample invertible sheaves on $X$ and $Y$ respectively.
If $\gs \times 1$ is not quasi-unipotent, then by \cite[Lemma~3.2]{Ke1}, there 
exists $r > 1$, $c > 0$, and an integral curve $C$ on $X \times Y$ such that
\begin{equation}\label{eq:Cartesian-growth}
 (((\gs \times 1)^*)^m (\pi_1^*\calH_X \otimes \pi_2^* \calH_Y).C) \geq c r^{m} 
 \qquad \text{for all }  m \geq 0.
\end{equation}
But
\[
((\gs \times 1)^*)^m (\pi_1^*\calH_X \otimes \pi_2^* \calH_Y) = 
\pi_1^*(\gs^*)^m\calH_X \otimes \pi_2^* \calH_Y.
\]
Since $\gs$ is quasi-unipotent, the intersection numbers of the right hand side with
any curve $C$ must be bounded by a polynomial. 
This contradicts \eqref{eq:Cartesian-growth}.
So $\gs \times 1$ must be quasi-unipotent.
Similarly, $1 \times \tau$ is quasi-unipotent.
Thus by Theorem~\ref{th:multi-main}, the pair 
$(\pi_1^* \LL)_{\gs \times 1}^{}, (\pi_2^* \calM)_{1 \times \tau}^{}$
is NC-ample and thus the ring of interest is noetherian.
\end{proof}

\section{Gel'fand-Kirillov dimension}

In this section we generalize the results of \cite[$\S$6]{Ke1}, showing
that a noetherian twisted multi-homogeneous coordinate ring has integer GK-dimension.
We first fix hypotheses on the ring $B$.

\begin{hypothesis}\label{hyp:ringB}
Let $X$ be a projective scheme with $s$ commuting NC-ample bimodules $\{ \LL_{(i, \gs_i)} \}$.
Assume that the commutation relations
of $\{ (\LL_i)_{\gs_i} \}$
are compatible in the sense of \eqref{eq:Bergman}. 
Let $B = B(X; \{ (\LL_i)_{\gs_i} \})$ and suppose that $B$ is right noetherian.
\end{hypothesis}

If $B$ is the twisted multi-homogenous coordinate ring associated to an NC-ample set
of invertible bimodules, then the vanishing of cohomology in Definition~\ref{def:NC-ample}
allows one to control the dimension of $B_{\ov{\imath}}$ for $\ov{\imath} \geq \ov{\imath}_0$ for some
$\ov{\imath}_0 \in \NN^s$. 
We are not guaranteed such control on the ``edges'' $\oplus_j B_{(0, \dots, j, \dots, 0)}$.
Thus, it will be easier to study the GK-dimension of the ideal
$B_{\geq \ov{\imath}_0}$ rather then the GK-dimension of $B$.

\begin{lemma}\label{lem:GK-ideals}
Let $B$ satisfy \eqref{hyp:ringB}, and let $\ov{\imath} \in \NN^s$.
Then $\GKdim B = \GKdim  (B_{\geq \ov{\imath}})_B$.
\end{lemma}
\begin{proof}
%We may replace $B_{(0, \dots, 0)}$ with $k$ so that the hypotheses of 
%\cite[Theorem~2.10]{G-multi-filtered} are satisfied. So 
%$\GKdim B = \max \{ \GKdim B_{\geq \ov{\imath}}, \GKdim B/B_{\geq \ov{\imath}} \}$.
%Thus we need only show $\GKdim B/B_{\geq \ov{\imath}} \leq \GKdim B_{\geq \ov{\imath}}$.
%
If $\ov{\jmath} \geq \ov{\imath}$, then $B_{\geq \ov{\jmath}} \subseteq B_{\geq \ov{\imath}}$ and
$\GKdim B_{\geq \ov{\jmath}} \leq \GKdim B_{\geq \ov{\imath}}$. So we may assume that 
$\ov{\imath}$ %= (a, a, \dots, a)$ for some $a \in \NN$ and that $a$
is sufficiently large so that $\LL_{(1,\gs_1)}^{j_1} \dots \LL_{(s,\gs_s)}^{j_s}$ is
generated by global sections for $(j_1, \dots, j_s) \geq \ov{\imath}$ by 
Proposition~\ref{prop:multi-tfae}\eqref{prop:multi-tfae2}.
So $B_{\ov{\jmath}} \subseteq B_{\ov{\jmath} + \ov{\imath}}$ for all $\ov{\jmath} \in \NN^s$.

We may grade $B$ by 
$\{ B_{\leq (n, n, \dots, n)} / B_{\leq (n-1, n-1, \dots, n-1)} \colon n \in \NN \}$
and grade $B_{\geq \ov{\imath}}$ by 
$\{ B_{\leq (n, n, \dots, n) + \ov{\imath}} / B_{\leq (n-1, n-1, \dots, n-1) + \ov{\imath}} 
\colon n \in \NN \}$.
Then
\[
\GKdim B = \varlimsup_n \log_n \dim B_{\leq (n, n, \dots, n)}
\leq \varlimsup_n \log_n \dim B_{\leq (n, n, \dots, n)+\ov{\imath}} = \GKdim B_{\geq \ov{\imath}}
\]
since $B_{\geq \ov{\imath}}$ is finitely generated \cite[Lemma~6.1]{KL}. The other inequality
$\GKdim B_{\geq \ov{\imath}} \leq \GKdim B$ is trivial.
\end{proof}

We will need to use multi-Veronese subrings and
also generalize a standard lemma for graded rings to the multi-graded case.

\begin{definition}
Let $B$ be a $k$-algebra, finitely multi-graded by $\NN^s$ (that is, each
multi-graded piece is finite dimensional). Then the subring
\[
B^{(n_1, \dots, n_s)} = \bigoplus_{(i_1, \dots, i_s) \in \NN^s} B_{(n_1 i_1, \dots, n_s i_s)}
\]
is a \emph{multi-Veronese subring} of $B$.
\end{definition}

\begin{lemma}\label{lem:multiveronese}
Let $B$ be a $k$-algebra, finitely multi-graded by $\NN^s$.
\begin{enumerate}
\item\label{lem:multiveronese1} 
If $B$ has ACC on multi-graded right ideals, then $B$ is right noetherian.
\item\label{lem:multiveronese2} 
If $B$ is right noetherian, then the multi-Veronese subring $A=B^{(n_1, \dots, n_s)}$
is right noetherian for any $(n_1, \dots, n_s) \in \NN^s$.
\end{enumerate}
\end{lemma}
\begin{proof}
Both claims are simple generalizations of the graded case.
For \eqref{lem:multiveronese1}, one may  
see that the conclusion is implicit in the proof that a right multi-filtered
ring is right noetherian if its associated multi-graded ring is right noetherian
\cite[Theorem~1.5]{G-multi-filtered}.
The proof of \eqref{lem:multiveronese2} is as in \cite[Proposition~5.10(1)]{AZ},
noting that if $I$ is a multi-graded ideal of $A$, then $I = IB \cap A$.
\end{proof}

Now we may replace $B$ with a multi-Veronese.

\begin{lemma}\label{lem:replace-multiveronse}
Let $B$ satisfy \eqref{hyp:ringB}, and let $\ov{n} \in \NN$.
Then $B^{(\ov{n})}$ satisfies Hypothesis~\ref{hyp:ringB}
and $\GKdim B = \GKdim B^{(\ov{n})}$.
\end{lemma}
\begin{proof}
Let $\ov{n} = (n_1, \dots, n_s)$ and $A= B^{(\ov{n})}$.
We may assume that all $n_i > 0$.
For the first claim, we have already seen in Lemma~\ref{lem:multiveronese} that
$A$ is right noetherian. The bimodules $\{ \LL_{(i, \gs_i)}^{n_i} \}$ commute
compatibly because their commutation relations are compositions
of the commutation relations for $\{ \LL_{(i, \gs_i)} \}$.
The bimodules $\{ \LL_{(i, \gs_i)}^{n_i} \}$ are also NC-ample by Theorem~\ref{th:multi-main}.

Now choose $\ov{m}=(m_1,\dots,m_s) \in \NN^s$ such that 
$\LL_{(1,\gs_1)}^{j_1} \dots \LL_{(s,\gs_s)}^{j_s}$ is
generated by global sections for $(j_1, \dots, j_s) \geq \ov{m}$ by 
Proposition~\ref{prop:multi-tfae}\eqref{prop:multi-tfae2}.
Then for $m_i \leq j_i < n_i+m_i$ there are short exact sequences
\[
0 \to \mathcal{K}_{(j_1, \dots, j_s)} \to B_{(j_1, \dots, j_s)} \otimes \OO_X
\to \LL_{(1,\gs_1)}^{j_1} \dots \LL_{(s,\gs_s)}^{j_s} \to 0.
\]
Then tensoring with $\LL_{(1,\gs_1)}^{n_1 a_1} \dots
\LL_{(s,\gs_s)}^{n_s a_s}$
and taking cohomology, we have
\begin{multline*}
B_{(j_1, \dots, j_s)} \otimes 
H^0(\LL_{(1,\gs_1)}^{n_1 a_1} \dots \LL_{(s,\gs_s)}^{n_s a_s})
\to B_{(j_1 + n_1 a_1, \dots, j_s + n_s a_s)} \\ \to 
H^1(\mathcal{K}_{(j_1, \dots, j_s)} \otimes 
\LL_{(1,\gs_1)}^{n_1 a_1} \dots \LL_{(s,\gs_s)}^{n_s a_s}).
\end{multline*}
For $(a_1, \dots, a_s)$ sufficiently large, the rightmost cohomology group vanishes.
So there exists $\ov{b}$ such that 
\[
B_{\geq \ov{b}} \subseteq \sum_{1 \leq i \leq s} \sum_{0 \leq j_i < n_i+m_i} B_{(j_1, \dots, j_s)}A.
\]
Hence $B_{\geq \ov{b}}$ is a finite $A$-module, so 
\[
\GKdim B = \GKdim (B_{\geq \ov{b}})_B = \GKdim (B_{\geq \ov{b}})_A \leq \GKdim A
\]
by Lemma~\ref{lem:GK-ideals} and \cite[Corollary~5.4]{KL}.
\end{proof}

We may now generalize \cite[Theorem~6.1]{Ke1} to the multi-homogeneous case.

\begin{theorem}\label{th:GK}
Let $B$ satisfy Hypothesis~\ref{hyp:ringB}.
Then $\GKdim B$ is an integer
and
\[
\dim X + 1 \leq \GKdim B \leq 
s\left((\ell+1) \dim X + 1\right),
\]
where $s$ is the number of commuting bimodules, $\rho=\rho(X)$ is the Picard number
of $X$, and $\ell = 2 \left\lfloor \frac{\rho - 1}{2} \right\rfloor$.
\end{theorem}
\begin{proof}
By Lemma~\ref{lem:replace-multiveronse}, we may replace $B$ with a multi-Veronese;
hence, replacing each $\gs_i$ with $\gs_i^{m_i}$ for some $m_i$, we may assume
each $\gs_i$ is unipotent. That is, up to numerical equivalence, $\gs_i^{-1} \equiv
I + N_i \in \GL_\rho(\ZZ)$. We know $N_i^{\ell+1} = 0$ for
all $i$ \cite[Lemma~6.12]{Ke1}. (We choose to use $\gs_i^{-1}$ since
we will use Cartier divisors and if $\LL \cong \OO_X(D)$, then $\LL^{\gs} \cong
\OO_X(\gs^{-1}D)$.)

Since the set of bimodules is NC-ample, we may again replace $B$ with a multi-Veronese
and assume that $H^q(X, \LL_{\ov{\gs}}^{\ov{n}}) = 0$ for all $q > 0, n_i > 0$
where $\ov{n} = (n_1,\dots,n_s)$. Thus 
$\dim H^0(X, \LL_{\ov{\gs}}^{\ov{n}}) = \chi(\LL_{\ov{\gs}}^{\ov{n}})$ for $n_i > 0$.
So by the Riemann-Roch Theorem \cite[p.~ 361, Example~18.3.6]{Fu},
\begin{equation} \label{eq:RR}
\dim H^0(X, \LL_{\ov{\gs}}^{\ov{n}})
 = \sum_{j=0}^{\dim X} \frac{1}{j!} 
 \int_X ( (\LL_{\ov{\gs}}^{\ov{n}})^{\bullet j}) \cap \tau_{X,j}(\OO_X).
\end{equation}
where $\bullet j$ denotes $j$th self-intersection and the $\tau_{X,j}(\OO_X)$
are constant $j$-cycles. By Lemma~\ref{lem:GK-ideals}, we may ignore 
$\dim H^0(X, \LL_{\ov{\gs}}^{\ov{n}})$ when some $n_i = 0$.

Let $D_i$ be a Cartier divisor such that $\LL_i \cong \OO_X(D_i)$.
The action of $\gs^{-n_i}$ on Cartier divisors modulo numerical equivalence is given by
\cite[(4.2), (4.3)]{Ke1}
\begin{align*}
\gs_i^{-n_i}  &\equiv \sum_{c=0}^\ell \binom{n_i}{c}N_i^c \\
\sum_{m=0}^{n_i-1} \gs_i^{-m} &\equiv \sum_{d=0}^\ell \binom{n_i}{d+1}N_i^d.
\end{align*}
So up to numerical equivalence,
\[
\LL_{\ov{\gs}}^{\ov{n}} \equiv 
\sum_{a = 1}^s \left[ \left( \prod_{b=1}^{a-1} \left( \sum_{c=0}^\ell 
\binom{n_b}{c}N_b^c \right) \right)\cdot \left( \sum_{d=0}^\ell \binom{n_a}{d+1}N_a^d D_a \right) \right].
\]
Thus $\dim H^0(X, \LL_{\ov{\gs}}^{\ov{n}})$ is a polynomial in $n_i$, $i = 1, \dots, s$,
with the degree of $n_i$ at most $(\ell+1)\dim X$, since one has at most a $(\dim X)$th 
self-intersection.

Now let $B_{\geq (1,1,\dots,1)}$ have the filtration given by assuming each $n_i \leq n$. 
Then
$f(n)= \dim B_{(1,1,\dots,1) \leq (n_1,\dots,n_s) \leq (n,\dots,n)}$ is a polynomial in $n$
of degree at most $s((\ell+1)\dim X + 1)$. This is because summing over each $i=1,\dots, s$ adds
$1$ to the degree of $n_i$. Then the degree $f(n)$ is maximized if
$\dim H^0(X, \LL_{\ov{\gs}}^{\ov{n}})$ has a term of the form 
$n_1^{(\ell+1)\dim X}\dots n_s^{(\ell+1)\dim X}$, since in this case, $(\ell+1)\dim X + 1$ is added
to itself $s$ times.

 Thus $\GKdim B$ is an
integer with the desired upper bound by Lemma~\ref{lem:GK-ideals}. 
Now by Lemma~\ref{lem:NC-ample-then-s-ample}, $B$ has a twisted homogeneous coordinate
ring $C$ as a subring. Now $\dim X + 1 \leq \GKdim C$ \cite[Theorem~7.17]{Ke-filters}, so
the lower bound on $\GKdim B$ holds.
\end{proof}

Examining \cite[Theorem~6.1]{Ke1}, 
\cite[Theorem~7.17]{Ke-filters} we see that these bounds on $\GKdim B$
are not optimal for the case $s=1$. However, the notational difficulties of 
repeating the arguments of \cite[$\S$6]{Ke1} for $s$ bimodules seem to outweigh
the benefits, given that exact results can be given in the following specific cases.

\begin{proposition}\label{prop:GK-Rees}
Let $\LL_\gs$ be a $\gs$-ample invertible bimodule on a projective scheme $X$. 
Let $B = B(X; \LL_\gs)$ be generated in degree one.
Then \[
 \GKdim B[It] = \GKdim B + 1 \]
where $I= B_{> 0}$ is the irrelevant ideal.
\end{proposition}
\begin{proof}
%By Theorem~\ref{th:GK}, $\GKdim B[It]$ is an integer. We have
%$\GKdim B \leq \GKdim B[It] \leq \GKdim B[t] = \GKdim B + 1$ \cite[Proposition~3.5]{KL}.
By Lemma~\ref{lem:replace-multiveronse}, we may replace $\gs$ with some $\gs^m$
and assume that $\gs$ is unipotent, $\dim B_m = \chi(\LL_\gs^m)$ for $m \geq 1$,
 and $\dim B_m \leq \dim B_{m+1}$
for $m \geq 0$. Let $f(m) = \dim B_m$. Then $\GKdim B = \deg f + 1$ \cite[(6.4)]{Ke1}.
Filter $B[It]$ by $(B[It])_{(i,j) \leq (n,n)}, n \in \NN$.
Now $\dim (B[It])_{(i,j)} = \dim B_{i+j}$, so 
\[
g(n)=\sum_{i=0}^n \sum_{j=0}^n
f(i+j) = \dim (B[It])_{(i,j) \leq (n,n)}.
\]
Since $f(m)$ is a numerical polynomial, $\deg g = \deg f + 2$.
So $\GKdim B[It] = \deg f + 2 = \GKdim B + 1$.
\end{proof}

For general $k$-algebras $R, S$, we have $\GKdim (R \otimes_k S) \leq \GKdim R + \GKdim S$
\cite[Lemma~3.10]{KL}.
However, for 
the tensor product of a twisted homogeneous coordinate ring and a general $k$-algebra, 
we have equality, as in the commutative case.

\begin{proposition}\label{prop:GK-tensor}
Let $\LL_\gs$ be $\gs$-ample on a projective scheme $X$, and let $B = B(X; \LL_\gs)$.
Let $S$ be any $k$-algebra.
Then \[ \GKdim (B \otimes S) = \GKdim B + \GKdim S. \]
\end{proposition}
\begin{proof}
There exists a Veronese subalgebra $B^{(n)}$ of $B$ such that
$f(m) = \dim B_m^{(n)}$ is a polynomial for $m >0$ and $\GKdim B = \GKdim B^{(n)}$
 \cite[(6.3)--(6.4)]{Ke1}. We may also assume that $B^{(n)}$ is generated
in degree one \cite[Theorem~7.17]{Ke-filters}.
Let $V = B_0 \oplus B_1^{(n)}$. Then $\dim V^m$ is a polynomial in $m$,
so $\GKdim (B \otimes S) = \GKdim B + \GKdim S$ \cite[Proposition~3.11]{KL}.
\end{proof}

\section*{Acknowledgements}
The author would like to thank Daniel Chan for suggesting this topic and his
continued interest in it. We also thank J.T.~Stafford and Michel Van~ den Bergh for their comments.

\bibliographystyle{amsplain}
%\bibliography{ample}

%\begin{comment}
\providecommand{\bysame}{\leavevmode\hbox to3em{\hrulefill}\thinspace}

%\end{comment}

\end{document}